\documentclass[11pt]{amsart}
\usepackage[margin=1in]{geometry}
\usepackage{url}
\usepackage{hyperref}
\hypersetup{colorlinks=true, citecolor=blue}

\usepackage[lite]{amsrefs}
\usepackage{amssymb}
\usepackage{mathrsfs}
\usepackage[shortlabels]{enumitem}
\usepackage[all,cmtip]{xy}

\newcommand{\isom}{\cong}

\newcommand{\C}{\mathbb{C}}

\newcommand{\E}{\mathcal{E}}

\newcommand{\I}{\mathcal{I}}
\renewcommand{\O}{\mathcal{O}}

\newcommand{\mult}{\mathrm{mult}}
\newcommand{\wt}{\widetilde}
\newcommand{\Bl}{\mathrm{Bl}}

\newcommand{\DB}{\underline{\Omega}} 

\numberwithin{equation}{section}

\theoremstyle{plain}
\newtheorem{thm}[equation]{Theorem}

\newtheorem{lem}[equation]{Lemma}

\newtheorem{conj}[equation]{Conjecture}

\theoremstyle{definition}

\theoremstyle{remark}

\begin{document}

\title[Upper bound on the multiplicity of rational and Du Bois singularities]{Upper bound on the multiplicity of rational and Du Bois singularities}
\author{Sung Gi Park}
\address{Department of Mathematics, Princeton University, Fine Hall, Washington Road, Princeton, NJ 08544, USA}
\address{Institute for Advanced Study, 1 Einstein Drive, Princeton, NJ 08540, USA}
\email{sp6631@princeton.edu \,\,\,\,\,\,sgpark@ias.edu}


\date{\today}

\begin{abstract}
This paper resolves a question of Huneke and Watanabe by proving a sharp upper bound for the multiplicity of Du Bois singularities: at a point of a $d$-dimensional variety with Du Bois singularities and embedding dimension $e$, the multiplicity is at most $\binom{e}{d}$. Additionally, the result recovers the previously known upper bound for the multiplicity of rational singularities.
\end{abstract}

\maketitle

\section{Introduction}
\label{sec:intro}

In an effort to prove Fujita's conjecture on the freeness of adjoint linear systems, Helmke \cite{Helmke97} obtained upper bounds for the multiplicity of log canonical centers of a log canonical pair $(X,D)$. More precisely, let $X$ be a smooth variety of dimension $n$ and $(X,D)$ be a log canonical pair. When $Z_d$ is a union of $d$-dimensional log canonical centers of $(X,D)$, then \cite{Helmke97}*{Corollary 4.6} states that the multiplicity at any point $x\in Z_d$ is bounded by $\binom{n}{d}$ (i.e. $\mult_x(Z_d)\le \binom{n}{d}$). Furthermore, if $Z_d$ is a minimal log canonical center and $d>0$, then $\mult_x(Z_d)\le \binom{n-1}{d-1}$. Helmke developed a general method to approach Fujita's conjecture, utilizing the aforementioned upper bounds on multiplicities, which led to the proof of Fujita's conjecture in low-dimensional cases.

The celebrated result of Koll\'ar-Kov\'acs \cite{KK10} states that any union of log canonical centers has Du Bois singularities, and Kawamata's subadjunction formula \cite{Kawamata98} implies that a minimal log canonical center has rational singularities. Subsequently, Helmke's multiplicity bound has been expected to be generalized to varieties with Du Bois or rational singularities. For a variety $X$ of dimension $d$ at a point $x\in X$ with embedding dimension $e=\dim m_x/m_x^2$, using commutative algebra methods, Huneke-Watanabe \cite{HW15}*{Theorem 3.1} proved the upper bound $\mult_x(X)\le \binom{e-1}{d-1}$ when $X$ has rational singularities at $x\in X$, and asked if the upper bound $\mult_x(X)\le \binom{e}{d}$ is true when $X$ has Du Bois singularities at $x\in X$ (see \cite{HW15}*{Remark 3.4}). Our main theorem gives an affirmative answer to this question and, by the same methods, also recovers the bound for rational singularities.

\begin{thm}
\label{thm: multiplicity bound for rational and Du Bois singularities}
Let $x\in X$ be a point in a variety with Du Bois (resp. rational) singularities. Denote $e:=\dim m_x/m_x^2$ and $d:=\dim X$ at $x$. Then we have
$$
\mult_xX\le \binom{e}{d} \quad \left(\mathrm{resp.} \binom{e-1}{d-1}\right).
$$
\end{thm}

In positive characteristic, Huneke–Watanabe \cite{HW15} proved the analogous bounds $\binom{e-1}{d-1}$ for F-rational singularities and $\binom{e}{d}$ for F-pure singularities. Note, however, that the bound for F-pure singularities does not imply that for Du Bois singularities in characteristic zero; see \cite{Schwede09} for precise connections between the two settings. Shibata \cite{Shibata17} proved the upper bound for a Cohen-Macaulay normal variety with Du Bois singularities, mainly following the proof of \cite{HW15}.

The proof of Theorem \ref{thm: multiplicity bound for rational and Du Bois singularities} is based on the multiplicity bound in $1$-dimensional case, which appears in Helmke's work, with an argument arising from Koszul complexes in a resolution of singularities. Throughout this paper, a variety is a reduced separated scheme of finite type over $\C$. The embedding dimension of a variety $X$ at $x$ is the dimension of the cotangent space $m_x/m_x^2$ at $x$, where $m_x$ is the maximal ideal sheaf of $x\in X$.

Motivated by Theorem \ref{thm: multiplicity bound for rational and Du Bois singularities} and Shibata \cite{Shibata17a}*{Proposition 3.9}, we propose the following multiplicity bound in terms of minimal log discrepancies.

\begin{conj}
Let $x\in X$ be a point on a variety with log canonical singularities. Denote $e:=\dim m_x/m_x^2$ and $d:=\dim X$ at $x$. If $\mathrm{mld}_x(X)>k-1$ for some integer $0\le k\le d$, then
$$
\mult_xX\le \binom{e-k}{d-k}.
$$
Moreover, equality holds when $k=d-1$.
\end{conj}

Here, $\mathrm{mld}_x(X)$ denotes the minimal log discrepancy of $X$ at $x$, i.e. the infimum of the log discrepancies $a(E;X)$ over all divisors $E$ with center $x$. Note that \cite{Shibata17a} proves the statement when $X$ is klt and $\mathrm{mld}_x(X)$ is replaced by the log canonical threshold of the maximal ideal at $x\in X$. Theorem \ref{thm: multiplicity bound for rational and Du Bois singularities} confirms the case $k=0$, and a minor refinement of its proof yields the case $k=1$. The equality statement is known in low dimensions: for surfaces with klt singularities by Artin \cite{Artin66}*{Corollary 6}, and for threefolds with terminal singularities by Kakimi \cite{Kakimi}*{Theorem 2.1} and by Takagi-Watanabe \cite{TW04}*{Proposition 3.10}.

\section{Preliminaries}

\subsection{Multiplicities of points on curves}

Let $X$ be a separated scheme of finite type over $\C$ and $x\in X$ be a closed point. For the multiplicity $\mult_x X$, it suffices to consider when $X$ has no embedded points through the following procedure. By the Noetherian property, there exists the largest ideal sheaf $\I\subset \O_X$ which is supported on a proper subvariety of $X$ not containing an irreducible component of $X$. Then the subscheme $X_0$ defined by $\O_X/\I$ satisfies $S_1$, or equivalently, there is no embedded point. It follows easily from the definition of the multiplicity that $\mult_xX=\mult_xX_0$. 

Let $C$ be a curve (a separated scheme of finite type over $\C$ with dimension $1$) that satisfies $S_1$. Note that the blow-up $\Bl_x C$ of a closed point $x\in C$ satisfies $S_1$. Indeed, the blow-up $\Bl_x C$ is the scheme-theoretic closure of $C\setminus \left\{x\right\}$ in $\Bl_x C$, so that $\Bl_x C$ has no embedded point if $C\setminus \left\{x\right\}$ has no embedded point. 

\begin{lem}[\cite{Helmke97}*{Lemma 4.2}]
\label{lem: multiplicity of curve}
Let $C$ be a curve with no embedded point. Let $\mu:\Bl_xC\to C$ be the blow-up of a closed point $x\in C$ and $E_x:=m_x\O_{\Bl_xC}$ be the exceptional divisor. Assume that for some $k\ge 0$, we have a natural inclusion
$$
\mu_*\O_{\Bl_xC}(-kE_x)\subset \O_C
$$
in $\mathcal R_C:=\oplus_\eta \O_{C,\eta}$, the sheaf of rational functions of $C$ where $\eta$ ranges over the generic points of $C$. Then
$$
\mult_x C\le \binom{e+k-1}{k}
$$
where $e=\dim m_x/m_x^2$ is the embedding dimension of $C$ at $x$.
\end{lem}

Note that $\mu:\Bl_xC\to C$ is an isomorphism over $C\setminus \left\{x\right\}$, which implies that both $\mu_*\O_{\Bl_xC}(-kE_x)$ and $\O_C$ are subsheaves of $\mathcal R_C$. We include the proof of Helmke \cite{Helmke97} for the reader's convenience.

\begin{proof}
Let $V=m_x/m_x^2$ be the cotangent space of $C$ at $x$. By choosing $e$ elements of $m_x$ whose residue classes form a basis of the vector space $V$, we obtain a surjection of rings
$$
\Phi:\bigoplus^k_{i=0}Sym^iV\to \O_C/m_x^{k+1}.
$$
Consider the following Cartesian square:
\begin{displaymath}
\xymatrix{
{I}\ar[r]\ar[d]& {\oplus^k_{i=0}Sym^iV}\ar[d]^{\Phi}\\
{\mu_*\O_{\Bl_xC}(-kE_x)/m_x^{k+1}}\ar[r]&{\O_C/m_x^{k+1}}
}
\end{displaymath}
Since $\Phi$ is surjective, the first column is surjective. The elements in $m_x^{k+1}$ are sections of $$\mu_*\O_{\Bl_xC}(-(k+1)E_x),$$ and thus we have a surjection
$$
\phi:I\to \mu_*\left(\O_{\Bl_xC}(-kE_x)/\O_{\Bl_xC}(-(k+1)E_x)\right)\isom \mu_*\O_{E_x}(-kE_x).
$$
Note that for any nonzero element $t\in V$, $tI$ is contained in the kernel of $\phi$. Therefore
$$
\mult_xC=h^0(E_x,\O_{E_x}(-kE_x))\le \dim I-\dim tI=\dim\left(\ker(t:I\to tI)\right)
$$
and $\ker(t:I\to tI)=Sym^kV$, from which we conclude that $\mult_xC\le \binom{e+k-1}{k}$.
\end{proof}
 
\subsection{Koszul complexes associated to sets of Cartier divisors}
\label{sec: Koszul complexes}

Let $\E$ be a vector bundle of rank $r$ on $X$, a separated scheme of finite type over $\C$. Given a morphism $\phi:\E\to \O_X$, there exists a Koszul complex $K_\bullet(\phi)$:
$$
0\to \wedge^r\E\to \dots\to\E\to \O_X\to 0.
$$
Locally, let $e_1,\dots, e_r$ be a generator of $\E$, then each map $d:\wedge^j\E\to\wedge^{j-1}\E$ for $1\le j\le r$ is described by
$$
d(e_{i_1}\wedge\cdots \wedge e_{i_j})=\sum_{k=1}^j(-1)^{k+1}\phi(e_{i_k})e_{i_1}\wedge\cdots \wedge\widehat{e_{i_k}}\wedge\cdots \wedge e_{i_j}
$$
where $e_{i_1}\wedge\cdots\wedge e_{i_j}$ is a generator of $\wedge^j\E$.

We consider the case when $\E=\oplus_{i=1}^r\O_X(-D_i)$ such that $D_i$ is an effective Cartier divisor, and $\phi:\E\to \O_X$ is the sum of natural inclusions $\O_X(-D_i)\to \O_X$. Note that
$$
\wedge^j\E\isom \bigoplus_{1\le i_1<\dots<i_j\le r} \O_X(-D_{i_1}-\cdots-D_{i_j}).
$$
Consider the complex
$$
\cdots\to\bigoplus_{1\le i_1<\dots<i_j\le r}\O_X(-D_{i_1}-\cdots-D_{i_j})\xrightarrow{d}\bigoplus_{1\le i_1<\dots<i_{j-1}\le r}\O_X(-D_{i_1}-\cdots-D_{i_{j-1}})\to\cdots,
$$
whose differential $d$ on each component $\O_X(-D_{i_1}-\cdots-D_{i_j})$ is given by the alternating sum $\sum_{k=1}^j(-1)^{k+1}\iota_k$ of natural inclusions
$$
\iota_k:\O_X(-D_{i_1}-\cdots-D_{i_j})\to \O_X(-D_{i_1}-\cdots-\widehat{D_{i_k}}-\cdots-D_{i_j}).
$$
Recall that $\O_X(-D_{i_1}-\cdots-D_{i_j})$ is the sub-$\O_X$-module of the sheaf of total quotient rings of $\O_X$ and $\iota_k$ is induced by this inclusion. Then this complex is the Koszul complex $K_\bullet(\phi)$. 

Indeed, let $f_1,\dots,f_r$ be the local generators of $D_1,\dots ,D_r$. Then $\O_X(-D_{i_1}-\cdots-D_{i_j})$ represents an ideal sheaf generated by $f_{i_1}\cdots f_{i_j}$. Treating the generator $f_{i_1}\cdots f_{i_j}$ as $e_{i_1}\wedge\cdots\wedge e_{i_j}$, we recover the local description of the Koszul complex $K_\bullet(\phi)$.

In the rest of the paper, we make use of this Koszul complex obtained from a set of Cartier divisors $D_1,\dots, D_r$. We follow the convention that for $j\ge 0$, the $-j$-th degree component of $K_\bullet(\phi)$ is 
$$
K_{-j}(\phi):=\bigoplus_{1\le i_1<\dots<i_j\le r} \O_X(-D_{i_1}-\cdots-D_{i_j}).
$$

\section{Proof of Theorem \ref{thm: multiplicity bound for rational and Du Bois singularities} for rational singularities}

Let $x\in X$ be a point in a variety with rational singularities, $f:\Bl_xX\to X$ be the blow-up at $x$, and $g:\wt X\to \Bl_xX$ be a resolution of singularities. Denote by $\mu:=f\circ g:\wt X\to X$, and let $\wt E_x:=g^*E_x$ be the pullback of the exceptional divisor of the blow-up $f$.

For a general hyperplane section $H$ of $X$ through $x$, its strict transform $\wt H$ in $\wt X$ is a base point free Cartier divisor such that $\mu^*H=\wt H+\wt E_x$. Since $X$ has rational singularities, we have $R\mu_*\O_{\wt X}=\O_X$, which induces the natural identity
\begin{equation}
\label{eqn: rational sing}
R\mu_*\O_{\wt X}(-\wt H-\wt E_x)=\O_X(-H)  
\end{equation}
from the projection formula.

Take general hyperplane sections $H_1,\dots, H_{d-1}$ through $x\in X$. Denote $C:= H_1\cap\cdots\cap H_{d-1}$. Then we have $\mult_xX=\mult_xC$. Consider the Koszul complex $K_\bullet(\phi)$ obtained from the morphism $\phi:\oplus_{i=1}^{d-1}\O_{\wt X}(-\wt H_i-\wt E_x)\to \O_{\wt X}$ and the Koszul complex $K_\bullet(\psi)$ obtained from the morphism $\psi:\oplus_{i=1}^{d-1}\O_{\wt X}(-\wt H_i)\to \O_{\wt X}$. Since $\wt H$ is base point free, $\wt H_1,\dots, \wt H_{d-1}$ is a regular sequence and their intersection $\wt C:=\wt H_1\cap\cdots\cap \wt H_{d-1}$ is a smooth curve which maps surjectively to $C$. Therefore, we have $K_\bullet(\psi)\isom\O_{\wt C}$.

Consider the twisted complex $K_\bullet(\psi)\otimes \O_{\wt X}(-(d-1)\wt E_x)$. Then for $0\le j\le d-1$, we have
$$
K_{-j}(\psi)\otimes \O_{\wt X}(-(d-1)\wt E_x)=\bigoplus_{1\le i_1<\dots<i_j\le d-1} \O_{\wt X}(-\wt H_{i_1}-\cdots-\wt H_{i_j}-(d-1)\wt E_x).
$$
Then we have a natural inclusion
$$
K_{-j}(\psi)\otimes \O_{\wt X}(-(d-1)\wt E_x)\to K_{-j}(\phi)
$$
induced by the natural inclusion
$$
\O_{\wt X}(-\wt H_{i_1}-\cdots-\wt H_{i_j}-(d-1)\wt E_x)\to \O_{\wt X}(-\wt H_{i_1}-\cdots-\wt H_{i_j}-j\wt E_x)
$$
of subsheaves of the sheaf of the total quotient ring of $\O_{\wt X}$ for each summand. From the description of Koszul complexes in Section \ref{sec: Koszul complexes}, this induces a morphism of complexes
\begin{equation}
\label{eqn: map of complexes for rational singularities}
K_\bullet(\psi)\otimes \O_{\wt X}(-(d-1)\wt E_x)\to K_\bullet(\phi).
\end{equation}

Recall that $K_\bullet(\psi)\otimes \O_{\wt X}(-(d-1)\wt E_x)\isom \O_{\wt C}(-(d-1)\wt E_x)$. Additionally, we have the natural isomorphism of complexes
$$
R\mu_*(K_\bullet(\phi))\isom K_\bullet(\oplus_{i=1}^{d-1}\O_{X}(-H_i)\to \O_{X})
$$
from \eqref{eqn: rational sing}. Therefore, $R^0\mu_*$ of \eqref{eqn: map of complexes for rational singularities} induces a morphism
$$
\mu_*\O_{\wt C}(-(d-1)\wt E_x)\to \O_C.
$$
Note that $\mu:\wt C\to C$ factors through the blow-up $\Bl_x C$, which induces a natural morphism
$$
\O_{\Bl_x C}(-(d-1)E_x)\to g_*\O_{\wt C}(-(d-1)\wt E_x).
$$
Consequently, we obtain a morphism
\begin{equation}
\label{eqn: final inclusion for rational sing}
f_*\O_{\Bl_x C}(-(d-1)E_x)\to \O_C.
\end{equation}
Here, $C$ has no embedded point since it is cut out by a regular sequence of hyperplane sections $H_1,\dots, H_{d-1}$ on a Cohen-Macaulay variety $X$.

Going back to \eqref{eqn: map of complexes for rational singularities}, we have a commuting diagram of morphisms of complexes
\begin{displaymath}
\xymatrix{
{\O_{\wt X}(-(d-1)\wt E_x)}\ar[r]\ar[d]& {\O_{\wt X}}\ar[d]\\
{K_\bullet(\psi)\otimes \O_{\wt X}(-(d-1)\wt E_x)}\ar[r]&{K_\bullet(\phi)}
}
\end{displaymath}
where the top row is obtained by taking $\bullet=0$ of the bottom row. Suppose a local section $s$ of $f_*\O_{\Bl_x C}(-(d-1)E_x)$ maps to a local section $s'$ of $\O_C$ via \eqref{eqn: final inclusion for rational sing}. Then $s'$ lifts to a section of $\O_{\wt X}$ which can be considered as a rational section of $\O_{\wt X}(-(d-1)\wt E_x)$. This rational section restricts back to the section $s$ of $\O_{\Bl_x C}(-(d-1)E_x)$, because \eqref{eqn: map of complexes for rational singularities} is a natural isomorphism over $X\setminus x$. Therefore, \eqref{eqn: final inclusion for rational sing} is a natural inclusion of subsheaves of the sheaf of rational functions of $C$.

The embedding dimension of $C$ at $x$ is $e-d+1$. Therefore, by Lemma \ref{lem: multiplicity of curve}, we obtain the multiplicity bound $\binom{e-1}{d-1}$ as desired.

\section{Proof of Theorem \ref{thm: multiplicity bound for rational and Du Bois singularities} for Du Bois singularities}

Let $x\in X$ be a point in a variety with Du Bois singularities. In the neighborhood of $x$, consider a closed embedding $X\subset W$ into a smooth variety $W$. Let $f:\Bl_xW\to W$ be the blow-up of $W$ at $x$, and $g:\wt W\to \Bl_xW$ be a log resolution of singularities such that the set-theoretic pullback $(f\circ g)^{-1}(X)=E+F$ is a reduced simple normal crossing divisor with $E=\mu^{-1}(x)$. Denote by $\mu:=f\circ g:\wt W\to W$, and let $\wt E_x:=g^*E_x$ be the pullback of the exceptional divisor of the blow-up $f$. Then $E$ is a reduced divisor of $\wt E_x$.

As in the proof for rational singularities, choose a general hyperplane section $H$ of $W$ through $x$. Then $\mu^*H=\wt H+\wt E_x$, where $\wt H$ is a base point free Cartier divisor. Since $W$ is smooth, we have an identity
\begin{equation}
\label{eqn: rational sing 2}
R\mu_*\O_{\wt W}(-\wt H-\wt E_x)=\O_W(-H).
\end{equation}
Additionally, $\mu$ is an isomorphism over the complement of $X$. Thus, we have a natural quasi-isomorphism of Du Bois complexes
$$
R\mu_*\DB_{\wt W,E\cup F}^0=\DB_{W,X}^0.
$$
See, for example, Du Bois \cite{DuBois81}*{Proposition 4.11}. Here, $\DB^0_{W,X}$ is the mapping cone of the natural morphism $\rho:\DB^0_W\to\DB^0_X$ of Du Bois complexes shifted by $[-1]$. In particular, if both $W$ and $X$ are Du Bois, then $\DB^0_{W,X}=\I_{W,X}$ where $\I_{W,X}$ is the ideal sheaf of $X\subset W$. This implies the identity
\begin{equation}
\label{eqn: Du Bois sing}
R\mu_*\O_{\wt W}(-F-E)=\I_{W,X}
\end{equation}
 (see also \cite{Park23}*{Corollary 1.10}).

Take general hyperplane sections $H_1,\dots, H_{d-1}$ through $x\in W$, and denote by $C:= X\cap H_1\cap\cdots\cap H_{d-1}$. Then we have $\mult_xX=\mult_xC$. Consider the Koszul complex $K_\bullet(\phi)$ obtained from the morphism 
$$
\phi:\O_{\wt W}(-F-E)\oplus\bigoplus_{i=1}^{d-1}\O_{\wt W}(-\wt H_i-\wt E_x)\to \O_{\wt W}
$$ 
and the Koszul complex $K_\bullet(\psi)$ obtained from the morphism
$$
\psi:\O_{\wt W}(-F)\oplus\bigoplus_{i=1}^{d-1}\O_{\wt W}(-\wt H_i)\to \O_{\wt W}.
$$
Note that $F,\wt H_1,\dots, \wt H_{d-1}$ is a regular sequence and their intersection $V:=F\cap\wt H_1\cap\cdots\cap\wt H_{d-1}$ is a reduced simple normal crossing divisor in a smooth variety $\wt H_1\cap\cdots\cap\wt H_{d-1}$. In particular, we have $K_\bullet(\psi)\isom\O_{V}$ and it is easy to see that $V$ is surjectively mapped to $C$.

Consequently, we have a morphism of complexes
\begin{equation}
\label{eqn: map of complexes for Du Bois singularities}
K_\bullet(\psi)\otimes \O_{\wt W}(-d\wt E_x)\to K_\bullet(\phi),
\end{equation}
induced by the natural inclusions
$$
\O_{\wt W}(-\wt H_{i_1}-\cdots-\wt H_{i_j}-d\wt E_x)\to \O_{\wt W}(-\wt H_{i_1}-\cdots-\wt H_{i_j}-j\wt E_x)
$$
and
$$
\O_{\wt W}(-F-\wt H_{i_1}-\cdots-\wt H_{i_j}-d\wt E_x)\to \O_{\wt W}(-F-\wt H_{i_1}-\cdots-\wt H_{i_j}-E-j\wt E_x)
$$
of subsheaves of the sheaf of the total quotient ring of $\O_{\wt W}$ for each $0\le j\le d$.

From \eqref{eqn: rational sing 2} and \eqref{eqn: Du Bois sing}, we have $R^0\mu_*(K_\bullet(\phi))=\O_{C}$. Recall that $K_\bullet(\psi)\otimes \O_{\wt W}(-d\wt E_x)\isom\O_{V}(-d\wt E_x)$. Therefore, $R^0\mu_*$ of \eqref{eqn: map of complexes for Du Bois singularities} induces a morphism
$$
\mu_*\O_{V}(-d\wt E_x)\to \O_C.
$$
Since $\mu:V\to C$ factors through the blow-up $\Bl_xC$, we have a morphism
$$
f_*\O_{\Bl_x C}(-dE_x)\to \O_C.
$$
Note that by the general choice of hyperplane sections, $C$ is reduced away from $x$. Let $C_{red}$ be the reduced scheme of $C$. Then we deduce a morphism
\begin{equation}
\label{eqn: final inclusion for Du Bois sing}
f_*\O_{\Bl_x C}(-dE_x)\to \O_{C_{red}}.
\end{equation}

Consider a commuting diagram of morphisms of complexes
\begin{displaymath}
\xymatrix{
{\O_{\wt W}(-d\wt E_x)}\ar[r]\ar[d]& {\O_{\wt W}}\ar[d]\\
{K_\bullet(\psi)\otimes \O_{\wt W}(-d\wt E_x)}\ar[r]&{K_\bullet(\phi)}
}
\end{displaymath}
where the top row is obtained by taking $\bullet=0$ of the bottom row. Suppose a local section $s$ of $f_*\O_{\Bl_x C}(-dE_x)$ maps to a local section $s'$ of $\O_{C_{red}}$ via \eqref{eqn: final inclusion for Du Bois sing}. Then $s'$ lifts to a section of $\O_{\wt W}$ which can be considered as a rational section of $\O_{\wt W}(-d\wt E_x)$. This rational section restricts to a rational section $s_V$ of $\O_{V}(-d\wt E_x)$. Because \eqref{eqn: map of complexes for Du Bois singularities} is a natural isomorphism over $W\setminus x$, the rational section $s_V$ is equal to the section $s$ pulled back to $V$ over the complement of $x$. Therefore, $s$ and $s'$ are the same over the complement of $x$, which indicates that \eqref{eqn: final inclusion for Du Bois sing} is a natural inclusion of subsheaves of the sheaf of rational functions of $C_{red}$.

The embedding dimension of $C$ at $x$ is $e-d+1$. Therefore, by Lemma \ref{lem: multiplicity of curve}, we obtain the multiplicity bound $\binom{e}{d}$ as desired.

\section{Boundary examples}

In this section, we provide an explicit example of a Du Bois singularity $x\in X$ of dimension $d$ and embedding dimension $e$ with $\mult_xX=\binom{e}{d}$. Let $\mathbb A^e$ be the affine space with coordinates $x_1,\dots,x_e$. Let $X\subset \mathbb A^e$ be the union of $d$-dimensional planes spanned by $d$ coordinate axes. Hence, $X$ is a union of $\binom{e}{d}$ planes of dimension $d$. It is well known and easy to see using inductive arguments that $X$ has Du Bois singularities, and $\mult_xX=\binom{e}{d}$ as desired. See \cite{HW15}*{Section 4} for a boundary example for rational singularities.


\begin{bibdiv}
    \begin{biblist}

\bib{Artin66}{article}{
   author={Artin, Michael},
   title={On isolated rational singularities of surfaces},
   journal={Amer. J. Math.},
   volume={88},
   date={1966},
   pages={129--136},
   issn={0002-9327},
   review={\MR{0199191}},
   doi={10.2307/2373050},
}

\bib{DuBois81}{article}{
   author={Du Bois, Philippe},
   title={Complexe de de Rham filtr\'{e} d'une vari\'{e}t\'{e} singuli\`ere},
   language={French},
   journal={Bull. Soc. Math. France},
   volume={109},
   date={1981},
   number={1},
   pages={41--81},
   issn={0037-9484},
   review={\MR{613848}},
}

\bib{Helmke97}{article}{
   author={Helmke, Stefan},
   title={On Fujita's conjecture},
   journal={Duke Math. J.},
   volume={88},
   date={1997},
   number={2},
   pages={201--216},
   issn={0012-7094},
   review={\MR{1455517}},
   doi={10.1215/S0012-7094-97-08807-4},
}

\bib{HW15}{article}{
   author={Huneke, Craig},
   author={Watanabe, Kei-ichi},
   title={Upper bound of multiplicity of F-pure rings},
   journal={Proc. Amer. Math. Soc.},
   volume={143},
   date={2015},
   number={12},
   pages={5021--5026},
   issn={0002-9939},
   review={\MR{3411123}},
   doi={10.1090/proc/12851},
}

\bib{Kakimi}{article}{
   author={Nobuyuki Kakimi},
   title={On the multiplicity of terminal singularities on threefolds},
   journal={preprint arXiv:math/0004105},
   date={2001},
}

\bib{Kawamata98}{article}{
   author={Kawamata, Yujiro},
   title={Subadjunction of log canonical divisors. II},
   journal={Amer. J. Math.},
   volume={120},
   date={1998},
   number={5},
   pages={893--899},
   issn={0002-9327},
   review={\MR{1646046}},
}

\bib{KK10}{article}{
   author={Koll\'{a}r, J\'{a}nos},
   author={Kov\'{a}cs, S\'{a}ndor J.},
   title={Log canonical singularities are Du Bois},
   journal={J. Amer. Math. Soc.},
   volume={23},
   date={2010},
   number={3},
   pages={791--813},
   issn={0894-0347},
   review={\MR{2629988}},
   doi={10.1090/S0894-0347-10-00663-6},
}

\bib{Park23}{article}{
   author={Park, Sung Gi},
   journal = {preprint arXiv:2311.15159v3},
   title = {{Du Bois complex and extension of forms beyond rational singularities}},
   year = {2023}
}

\bib{Shibata17}{article}{
   author={Shibata, Kohsuke},
   title={Upper bound of the multiplicity of a Du Bois singularity},
   journal={Proc. Amer. Math. Soc.},
   volume={145},
   date={2017},
   number={3},
   pages={1053--1059},
   issn={0002-9939},
   review={\MR{3589305}},
   doi={10.1090/proc/13307},
}

\bib{Shibata17a}{article}{
   author={Shibata, Kohsuke},
   title={Multiplicity and invariants in birational geometry},
   journal={J. Algebra},
   volume={476},
   date={2017},
   pages={161--185},
   issn={0021-8693},
   review={\MR{3608148}},
   doi={10.1016/j.jalgebra.2016.11.027},
}

\bib{Schwede09}{article}{
   author={Schwede, Karl},
   title={$F$-injective singularities are Du Bois},
   journal={Amer. J. Math.},
   volume={131},
   date={2009},
   number={2},
   pages={445--473},
   issn={0002-9327},
   review={\MR{2503989}},
   doi={10.1353/ajm.0.0049},
}

\bib{TW04}{article}{
   author={Takagi, Shunsuke},
   author={Watanabe, Kei-ichi},
   title={On F-pure thresholds},
   journal={J. Algebra},
   volume={282},
   date={2004},
   number={1},
   pages={278--297},
   issn={0021-8693},
   review={\MR{2097584}},
   doi={10.1016/j.jalgebra.2004.07.011},
}


    \end{biblist}
\end{bibdiv}
\end{document}